\documentclass[letterpaper, 10 pt, conference]{ieeeconf}  

\IEEEoverridecommandlockouts                              

\usepackage{amsmath}
\usepackage{color}
\usepackage{amssymb}  
\usepackage{amsfonts}
\usepackage{graphicx}
\usepackage{dsfont}
\usepackage{psfrag}
\usepackage{epstopdf}
\usepackage{mathtools}
\usepackage[noadjust]{cite}
\usepackage{bbm}
\usepackage{algorithm}
\usepackage{algorithmic}
\usepackage{savesym}
\usepackage{amsthm}

\theoremstyle{definition}

\newtheorem{remark}{Remark}

\theoremstyle{plain}

\title{\LARGE \bf
Fractional-Order Model Predictive Control for Neurophysiological Cyber-Physical Systems: A Case Study using Transcranial Magnetic Stimulation
}

\author{Orlando Romero$^{\dag *}$ \qquad Sarthak Chatterjee$^{\ddag *}$ \qquad S\'ergio Pequito$^{\dag}$
\thanks{$^{\dag}$Orlando Romero and S\'ergio Pequito are with the Department of Industrial and Systems Engineering, Rensselaer Polytechnic Institute, Troy, NY 12180, USA
        {\tt\small \{rodrio2,goncas\}@rpi.edu}}%
\thanks{$^{\ddag}$Sarthak Chatterjee is with the department of Electrical, Computer, and Systems Engineering, Rensselaer Polytechnic Institute, Troy, NY 12180, USA
        {\tt\small chatts3@rpi.edu}}%
\thanks{$^*$Equal contribution}        
        }

\begin{document}

\maketitle
\thispagestyle{empty}
\pagestyle{empty}

\begin{abstract}

Fractional-order dynamical systems are used to describe processes that exhibit temporal long-term memory and power-law dependence of trajectories. There has been evidence that complex neurophysiological signals like electroencephalogram (EEG) can be modeled by fractional-order systems. In this work, we propose a model-based approach for \mbox{closed-loop} Transcranial Magnetic Stimulation (TMS) to regulate brain activity through EEG data. More precisely, we propose a model predictive control (MPC) approach with an underlying fractional-order system (FOS) predictive model. Furthermore, MPC offers, by design, an additional layer of robustness to compensate for system-model mismatch, which the more traditional strategies lack. To establish the potential of our framework, we focus on epileptic seizure mitigation by computational simulation of our proposed strategy upon seizure-like events. We conclude by empirically analyzing the effectiveness of our method, and compare it with event-triggered open-loop strategies.

\end{abstract}

\section{Introduction}

In the context of neurophysiological signals, \emph{temporal} fractional properties in both health and disease states have become apparent and, with it, fractional-order systems have demonstrated a huge potential for clinical \mbox{applications~\cite{brodu2012exploring,ihlen2010interaction,ciuciu2012scale,zorick2013multifractal,zhang2015multifractal,papo2014functional,suckling2008endogenous}}. Practically, this leads signals to become non-stationary and to possess long-term memory dependencies with themselves, with the backwards-decaying weights of such dependencies following a power-law distribution~\cite{1_tarasov2011fractional}. The persisting temporal dependencies illustrated by these systems have given rise to \emph{fractional-order} based modeling, design, and analysis of novel neurotechnologies.

Recently, dynamical \emph{spatiotemporal} fractional models have been proposed as a tool to model neurophysiological signals suitable to deal with structured data and to equip us with modeling capabilities that capture spatial (i.e., the contributions of the signal's components into each other) and temporal long-range memory through the \mbox{so-called} fractional-order coefficients associated with the power-law exponents~\cite{13_xue2016spatio,14_xue2016minimum,xue2017reliable,15_xue2017constructing,PequitoC24,pequitoC29,jornalVassilisYuankun}, and possibly under unknown unknowns~\cite{gaurav2017acc,gupta2018re}.

Notwithstanding the above, the main advent of \mbox{model-based} approaches is that we can understand how an external signal or stimulus would craft the dynamics of the process. Simply speaking, it enables us to \emph{design} a sequence of interactions (i.e., a control strategy or law) with the system such that we can steer its dynamics towards satisfying desirable properties. That said, due to the highly dynamic nature of neurophysiological processes, it is imperative that we consider \emph{feedback} mechanisms~\cite{andrzej}. In other words, we need to leverage the continuous flow of measurements of the system to tune (for the individual's process) the control strategy. A particularly successful strategy that has achieved remarkable success in several engineering applications is the strategy of \emph{model predictive control} (MPC) that consists of three key ideas~\cite{camacho2013model,mayne2000constrained,qin2003survey,allgower2012nonlinear,garcia1989model,morari1999model}: (\emph{i}) a model-based approach; (\emph{ii}) capability of predicting the evolution of the system and its states upon a devised feedback control strategy that aims to optimize an objective that encapsulates the risk assessment of abnormal behavior; and (\emph{iii}) receding finite-horizon re-evaluation of the control strategy performance devised in the previous point.

In the context of neurophysiological processes, we propose to leverage fractional-order models to equip us with the aforementioned prediction and control capabilities that go hand in hand with the closed-loop design of Cyber-Physical Systems (CPSs). As a consequence, we will be able to develop stimulation strategies in the form of Transcranial Magnetic Stimulation (TMS) from electroencephalographic (EEG) data, and use it to annihilate or mitigate the overall duration and strength of an epileptic seizure. Notwithstanding, we believe that similar design and strategies can be envisioned in other contexts where closed-loop \mbox{deep-brain} electrical stimulation is available, e.g., Parkinson's disease~\cite{schuepbach2013neurostimulation,witt2013relation}, Alzheimer's disease~\cite{nardone2015neurostimulation}, depression~\cite{marangell2007neurostimulation,bewernick2010nucleus}, and anxiety~\cite{sturm2007nucleus}, just to mention a few.

To summarize, in what follows we introduce in a pedagogical manner the control mechanisms to be deployed as part of future neurophysiological cyber-physical systems, with particular emphasis on TMS for epilepsy. Ultimately, the integration of these design features will lead to more reliable CPSs that will immediately improve, or otherwise bring a positive impact, on the quality of life of the patients that qualify for the use of such technologies.

The paper is organized as follows. Section~\ref{sec:prelim} introduces the preliminary setup of our problem. Section~\ref{sec:fos-mpc} presents the fractional-order model predictive control framework. Finally, we present illustrative examples demonstrating the efficacy of epileptic seizure mitigation using TMS in Section~\ref{sec:simulations}.

\section{Preliminaries}\label{sec:prelim}


There are primarily two kinds of brain stimulation strategies that exist, which include \emph{open-loop} and \emph{closed-loop} strategies. Open-loop stimulation strategies consist of any stimulation strategy that does not utilize current brain activity data to regulate the stimuli applied to the patient's brain. On the other hand, closed-loop (brain-responsive) strategies consist of stimulation treatments based on automatic electrical stimulation directly influenced by the present (i.e., \mbox{real-time}) behavior being observed through continuous recording of brain activity. Past data can also be used in this strategy, but due to limited storage capabilities, these mechanisms are usually designed to depend exclusively on the actual recorded and stored data at a given time, consisting of a finite temporal window ranging from a fixed number of past instances of discretized time to the present measurement.

\subsection{Transcranial Magnetic Stimulation}

Transcranial Magnetic Stimulation (TMS) is a noninvasive form of brain stimulation created by inducing electric currents at specific areas of the brain by the help of a changing magnetic field using electromagnetic induction. This is done using an electric pulse generator connected to a magnetic coil, which in turn, is connected to the scalp. Although research on this form of neurostimulation is still evolving, it has been shown to demonstrate therapeutic potential in neurodegenerative disorders like Alzheimer's disease~\cite{lefaucheur2014evidence}, motor neuron disease~\cite{fang2013repetitive}, and stroke~\cite{dimyan2010contribution}.

\subsection{Linear fractional-order systems} 
For many biological systems, linear time-invariant (LTI) \mbox{state-space} models are insufficient to accurately capture the real evolution of the systems for anything other than a very small interval of time into the future, given that the current state of the system may have a non-negligible dependence on several past states, or even from the states ranging from the entire period of time so far. In this paper, we focus solely on discrete-time systems and control, and we will model the above scenario as
\begin{equation}
    x_{k+1} = f_k(x_k,x_{k-1},\ldots,x_0) + w_k,
\end{equation}
with known functions $f_k:\mathbb{R}^{n(k+1)}\to\mathbb{R}^n$ for $k\in\mathbb{Z}_+$ and $w_k$ being the \emph{process noise}. Alternatively, we may consider finite-history models of the form 
\begin{equation}
    x_{k+1} = f_k(x_k,\ldots,x_{k-p+1}) + w_k,
\end{equation}
such as multivariate autoregressive (MVAR) models
\begin{equation}
    x_{k+1} = \sum_{j=0}^{p-1}A_jx_{k-j} + w_k,
    \label{eq:MVAR}
\end{equation}
with $A_0,\ldots,A_{p-1}\in\mathbb{R}^{n\times n}$. However, system identification for such autoregressive models can become mathematically intractable or suffer issues of numerical instability that originate from the possibly large number of parameters to be estimated ($n^2p$ entries in $A_0,\ldots,A_{p-1}$, plus usually a certain other number for the covariance matrix of the noise).

For the reasons covered above, we introduce the so-called linear \emph{fractional-order system} (FOS) models, of the form
\begin{equation}
    \Delta^\alpha x_{k+1} = Ax_k + w_k,
\end{equation}
where the \emph{state coupling matrix} $A\in\mathbb{R}^{n\times n}$ and the vector of \emph{fractional-order} exponents $\alpha\in\mathbb{R}_+^n$ are now the only parameters to be estimated. This model can be rewritten~\cite{guermah2008controllability} as 
\begin{equation}
    x_{k+1} = \sum_{j=0}^k A_jx_{k-j} + w_k,
\end{equation}
which can be readily approximated as an MVAR model akin to~\eqref{eq:MVAR}. Furthermore, for the simulations conducted in this paper, we will model $w_k$ as an additive white Gaussian noise (AWGN).

\section{Fractional-Order Model Predictive Control} \label{sec:fos-mpc}

We start by showing that~\eqref{eq:MVAR} can be written as an \emph{augmented} LTI system model. To do this, let
\begin{equation}
    \tilde{x}_k = \begin{bmatrix}x_k\\ \vdots \\ x_{k-p+1}\end{bmatrix}
\end{equation}
denote the so-called \emph{augmented} state vector, with the understanding that $x_k=0$ for $k\in\{-1,\ldots,-p+1\}$. Then, clearly, the first block in $\tilde{x}_{k+1}$ can be expressed as a linear combination of the block in $\tilde{x}_k$ through \eqref{eq:MVAR}. On the other hand, the remaining $p-1$ blocks in $\tilde{x}_{k+1}$ precisely match the first $p-1$ blocks of $\tilde{x}_k$. More precisely, we have
\begin{align}
    \tilde{x}_{k+1} &= \underbrace{\begin{bmatrix} A_0 & \ldots & A_{p-2} & A_{p-1}\\ I & \ldots & 0 & 0\\ \vdots & \ddots & \vdots & \vdots\\ 0 & \ldots & I & 0\end{bmatrix}}_{=\tilde{A}}\tilde{x}_k + \underbrace{\begin{bmatrix}I \\ 0\\ \vdots \\ 0\end{bmatrix}}_{=\tilde{B}^w}w_k \nonumber \\&= \tilde{A}\tilde{x}_k + \tilde{B}^w w_k
\end{align}
for $k=0,1,\ldots$, which is an LTI model that we will refer to as the $p$-\emph{augmented} LTI system. 

It should also be noted that, if we consider a \emph{time-varying} FOS system
\begin{equation}
    \Delta^\alpha x_{k+1} = A_k x_k + w_k,
\end{equation}
then \eqref{eq:MVAR} would still be valid, except that $A_j$ would need to be indexed by $k$ as well and we could use the same reasoning to derive a finite-history approximation
\begin{equation}
    \tilde{x}_{k+1} = \tilde{A}_k\tilde{x}_k + \tilde{B}^w w_k,
\end{equation}
with the only exception that this $p$-augmented model is linear and \emph{time-varying}.

These representations will play a key role in the implementation of our proposed model predictive control strategy, since such an MVAR approximation of an FOS plant, \mbox{re-written} as an LTI model, allows for efficient numerical solutions to be determined for quadratic-cost optimal control. The reason for this is that a full representation of an FOS model would require increasing memory storage, while making the problem of computing optimal control actions an increasingly computationally demanding task. Naturally, such a completely faithful representation of \mbox{fractional-order} systems would therefore become intractable for most practical applications at likely little cost benefit, since the temporal dependence, while a long-term one for FOS models, is also decaying with respect to the weights.

\subsection{Model Predictive Control}
We now focus on the design of a full-state feedback controller for a linear time-varying system over discrete-time,
\begin{equation}
    x_{k+1} = A_kx_k + B_ku_k + B_k^ww_k,
\end{equation}
where $w_k$ denotes a sequence of independent and identically distributed (i.i.d.) random vectors, following a standard normal distribution with zero mean and identity covariance matrix. The weight matrices $B_k^w$ are intended to make the notation cleaner and add flexibility to the process noise $w_k' = B_k^w w_k$ since the covariance matrix of $w_k'$ does not need to be positive definite, which is a strict requirement for Gaussian distributions. The objective is to design the feedback controller such that it minimizes a quadratic cost functional of the input and state vectors over an infinite time horizon. In other words, the objective is to determine the sequence of control inputs $u_0,\ldots,u_{N-1}$ that minimizes a quadratic cost function of the form

\begin{equation}
\begin{aligned}
& \underset{{u}_0,\ldots,{u}_{N-1}}{\text{minimize}}
& & \mathbb{E}\left\{\sum_{k=1}^{N}\|{x}_k\|_{{Q}_k}^2 + \sum_{k=1}^N{c}_k^\mathsf{T}{x}_k + \sum_{k=0}^{N-1}\|{u}_k\|_{{R}_k}^2\right \} \\
& \text{subject to}
& & {x}_0 = {x}_\mathrm{init}, \\
& & & {x}_{k+1} = {A}_k{x}_k + {B}_k{u}_k + {B}_k^w{w}_k,\\
& & & \text{for} \:\: k=0,1,\ldots,N-1,\\
\end{aligned}
\label{eq:LQG}
\end{equation}
for $N\to\infty$, where $x_k \in \mathbb{R}^n$, $u_k \in \mathbb{R}^{n_u}$, and \mbox{${Q}_1,\ldots,{Q}_N\in\mathbb{R}^{n\times n}$} and ${R}_0,\ldots,{R}_{N-1}\in\mathbb{R}^{n_u\times n_u}$ are given positive semidefinite matrices. Recall that, $Q\in\mathbb{R}^{n\times n}$ is a \emph{positive semidefinite} matrix if $x^\mathsf{T} Qx\geq 0$ for every \mbox{$x\in\mathbb{R}^n$}, and $\|x\|_Q = \sqrt{x^\mathsf{T} Qx}$ in that case. For seizure mitigation via TMS, we propose to use $Q_k = I_{n\times n}$, \mbox{$c = 0_{n\times 1}$}, and $R_k = \varepsilon I_{n_u\times n_u}$ with $\varepsilon >0$, such that the objective becomes largely to steer the total energy in the expected value of the brain signals towards the smallest amount possible. The quadratic term on the input, which represents the stimulation signal induced in the brain via TMS, is intended to add a penalization term for stimulating the patient too harshly, since this may be unsafe, create discomfort for the patient, or have harmful psychological effects~\cite{moratti2014adverse}.

Problem~\eqref{eq:LQG} is often referred to as the \mbox{\emph{linear-quadratic-Gaussian} (LQG)} control problem~\cite{LOC}, and it may be solved offline and implemented in real time through a feedback control strategy of the form 
\begin{equation}
    {u}_k = {K}_k{x}_k\qquad (k=0,1,\ldots,N-1),
\end{equation}
for some appropriate matrices ${K}_0,\ldots,{K}_{N-1}\in\mathbb{R}^{n_u\times n}$ (called \emph{feedback gain} matrices)~\cite{Hespanha09}, which are functions of $x_0$. Furthermore, \eqref{eq:LQG} can be written as an unconstrained \emph{quadratic program} (QP), i.e., the minimization of a function $f(x) = \frac{1}{2}\|x\|_Q^2 + c^\mathsf{T} x$. Such problems can be efficiently solved numerically, for instance, via the \texttt{quadprog()} function in \textsc{Matlab}. In fact, even if we include additional \emph{linear} constraints on the state and input vectors, the problem will remain a QP. In particular, we may consider constraints of the form $u_\mathrm{min} \leq (u_k)_i \leq u_\mathrm{max}$ for every time step $k$ and every input component indexed by $i$. This can be used to efficiently ensure that our proposed TMS approach will only ever administer safe voltage stimuli. 

In model predictive control, the objective is to solve these consecutive constrained finite-horizon LQG problems over a moving horizon, in order to create additional robustness, as compared to directly solving for $N\to\infty$. More precisely, at time step $k$, the proposed strategy involves solving the problem

\begin{equation}
\begin{aligned}
& \underset{{u}_k,\ldots,{u}_{k+P-1}}{\text{minimize}}
& & \mathbb{E} \Bigg\{ \sum_{j=1}^{P}\|{x}_{k+j}\|_{{Q}_{k+j}}^2 + \sum_{j=1}^P{c}_{k+j}^\mathsf{T}{x}_{k+j} \\&&&+ \sum_{j=0}^{P-1}\|{u}_{k+j}\|_{{R}_{k+j}}^2 \Bigg\} \\
& \text{subject to}
& & {x}_k = \text{ observed or estimated current state}, \\
& & & {x}_{k+j+1} = {A}_{k+j}{x}_{k+j} + {B}_{k+j}{u}_{k+j} \\
& & & + {B}_{k+j}^w{w}_{k+j},\\
& & & \text{for} \:\: j=0,1,\ldots,P-1,\\
& & & \text{other linear constraints on } x_{k+1},\ldots,x_{k+P},\\&&&u_k,\ldots,u_{k+P-1},
\end{aligned}
\label{eq:mpc123}
\end{equation}
where $P$ is called the \emph{prediction horizon}, but it only deploys the control strategy associated with the first $M$ time steps (referred to as the \emph{control horizon}). Simply speaking, after we reach the state $x_{k+M-1}$, we update $k \leftarrow k+M-1$ and recompute the new LQG solution. This way, by design, there is a new layer of robustness that the solution of an LQG (even for $N\to\infty$) does not offer by itself, since the optimal strategy is constantly being re-evaluated based on short-term control action implementation of a long-term prediction~\cite{Bequette2013}.

\section{Simulation Results}\label{sec:simulations}

In what follows, we propose to illustrate the use of the fractional-order model predictive control framework for TMS in the context of annihilating or mitigating epileptic seizures using real EEG data. We start by considering an epileptic seizure captured by a linear FOS model whose parameters are obtained through a system identification method using real brainwave data. Then, we simulate three different stimulation strategies: \emph{(i)} an \mbox{open-loop} stimulation strategy; \emph{(ii)} an \mbox{event-triggered} open-loop strategy (i.e., once a certain event, in this case the beginning of a seizure, is detected, a \mbox{short-term} open-loop stimulation strategy is deployed); and \emph{(iii)} a closed-loop stimulation strategy based on our proposed \mbox{FOS-MPC} approach. 

\subsection*{System Identification Through Parameter Estimation with EEG data}
First, we need to determine the parameters $A$ and $\alpha$ that model both spatial coupling and fractional exponents, respectively, that craft the evolution of the state $x_k\in\mathbb{R}^n$ in the FOS model
\begin{equation}
    \Delta^\alpha x_{k+1} = Ax_k + w_k,
    \label{eq:sysIDed}
\end{equation}
with $w_k$ denoting additive white Gaussian noise (AWGN). We will consider only $n=4$ components in our state vector $x_k$ with each one corresponding to the measurement obtained from a single EEG microelectrode. To identify the parameters $A$, $\alpha$, and the covariance matrix $\Sigma_w$ (assumed equal to the diagonal matrix $\Sigma = \sigma_w^2 I_{4\times 4}$ for simplicity) associated with $w_k$, we used the method proposed in~\cite{gaurav2017acc}. For illustration purposes, we consider the normalized (i.e., voltage units have been scaled to the interval $-1$ to $1$) EEG recordings 1-4 of subject 11 from the CHB-MIT Scalp EEG database~\cite{PhysioNet}. The parameters obtained are as follows:
\begin{subequations}\label{eq:concreteparam}
\begin{equation}
    A = \begin{bmatrix}    
    0.0402  &  0.0604  & -0.0040  & -0.0450\\
    0.0340  & -0.0571  &  0.0742  &  0.0701\\
   -0.0119  & -0.0032  & -0.0105  &  0.0078\\
   -0.0335  &  0.0165  & -0.0009  &  0.0453\end{bmatrix},
\end{equation}
\begin{equation}
    \alpha = \begin{bmatrix}0.6606  &  0.7973  &  1.0670  &  0.6977\end{bmatrix}^\mathsf{T},
\end{equation}
\end{subequations}
and $\sigma_w^2 = 0.2$. More precisely, we utilized a normalized \mbox{10-second} sample during a period of ictal activity (i.e., activity during an epileptic seizure), sampled at 160 Hz.

\subsection*{Experiment 1: Open-Loop Electrical Neurostimulation}
Based upon the identified parameters of the \mbox{system~\eqref{eq:sysIDed}--\eqref{eq:concreteparam}}, we consider the following forced FOS 
\begin{equation}
    \Delta^\alpha x_{k+1} = Ax_k + Bu_k + w_k + d_k,
\end{equation}
where $d_k$ denotes a disturbance term that can be understood as a persistent neural activity incoming from nearby regions. For the sake of our simulations, $d_k$ will be set as frequent wavelet-like bursts of amplitude $d_k=0.25$ and $d_k=1$, which occur in a disjoint manner at random points of time according to Poisson counting processes with rates of 0.2 s and 1 s, respectively.

\begin{figure}[!ht]
    \centering
    \includegraphics[width=0.5\textwidth]{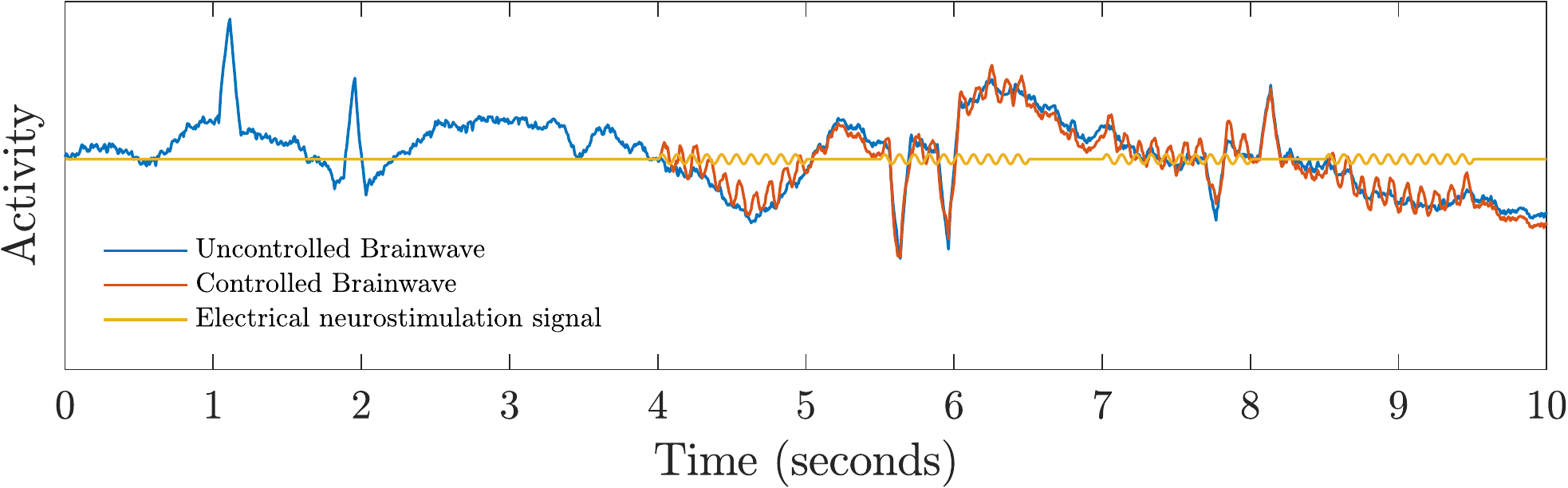}
    \caption{Simulated seizure is depicted (in blue) versus an attempted controlled signal (in red) under an open-loop stimulation strategy (in yellow).}
    \label{fig:openloop}
\end{figure}

Additionally, $u_k$ denotes the voltage stimulus signal being applied. For the sake of simplicity, we will assume \mbox{$B = [1,1,1,1]^\mathsf{T}$}, which corresponds to a stimulus that perturbs all the channels uniformly. 

Several open-loop stimulation strategies were implemented, which always entail $u_k$ to be determined without information about the state of the system. Specifically, when considering the previous model, we start actuation at the 4 second mark, when the EEG signal (blue) branches out into two signals, one being the continuation of the original signal without any stimulation, and the other being the electrically stimulated version (red) of the signal through the open-loop controlled stimulus (orange). Here, the input will consist of 1 s periods of sinusoidal activity of amplitude 0.5, frequency 16 Hz, phase 0, and consecutively followed by unstimulated periods lasting 0.5 s.

\begin{figure}[!ht]
    \centering
    \includegraphics[width=0.5\textwidth]{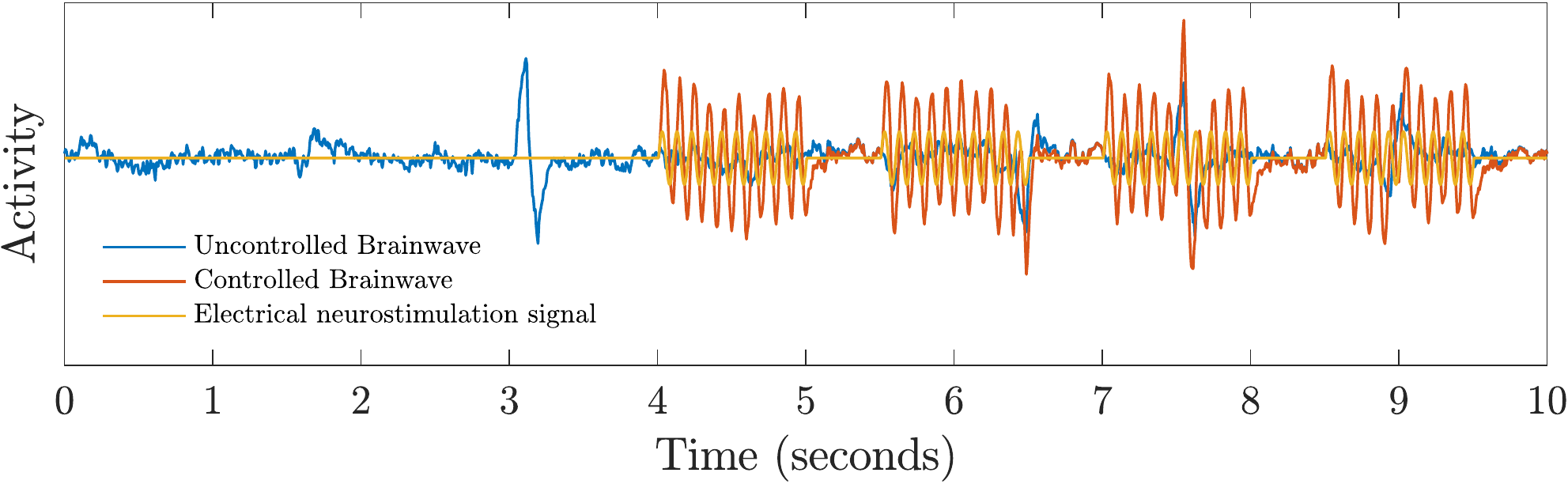}
    \caption{Seizure-like activity (in red)  is induced by an open-loop stimulation strategy (in yellow) of an otherwise regular activity depicted in blue.}
    \label{fig:openloopseizure}
\end{figure}

As we can see in Fig.~\ref{fig:openloop}, this strategy is unable to steer the evolution of the system towards normal activity. Specifically,  the controlled brainwaves seem to be almost unaffected, which is consistently observed if we adopt slightly different open-loop strategies (e.g., adopting a different waveform, as it is in the case of a biphasic rectangular pulse). In fact, we have observed that such open-loop control strategies can induce the increase of activity and, at times, even originate \mbox{seizure-like} activity. As a particular instance of these findings, we re-identified the parameters in our system based on inter-ictal data of the same subject as before. Next, we considered a stimulation strategy as depicted in Fig.~\ref{fig:openloopseizure}, which lead to several periods of seizure-like activity which would not have occurred if the original system was left to evolve on its own.

\subsection*{Experiment 2: Event-Triggered Open-Loop Control} 
\begin{figure}[!ht]
    \centering
    \includegraphics[width=0.5\textwidth]{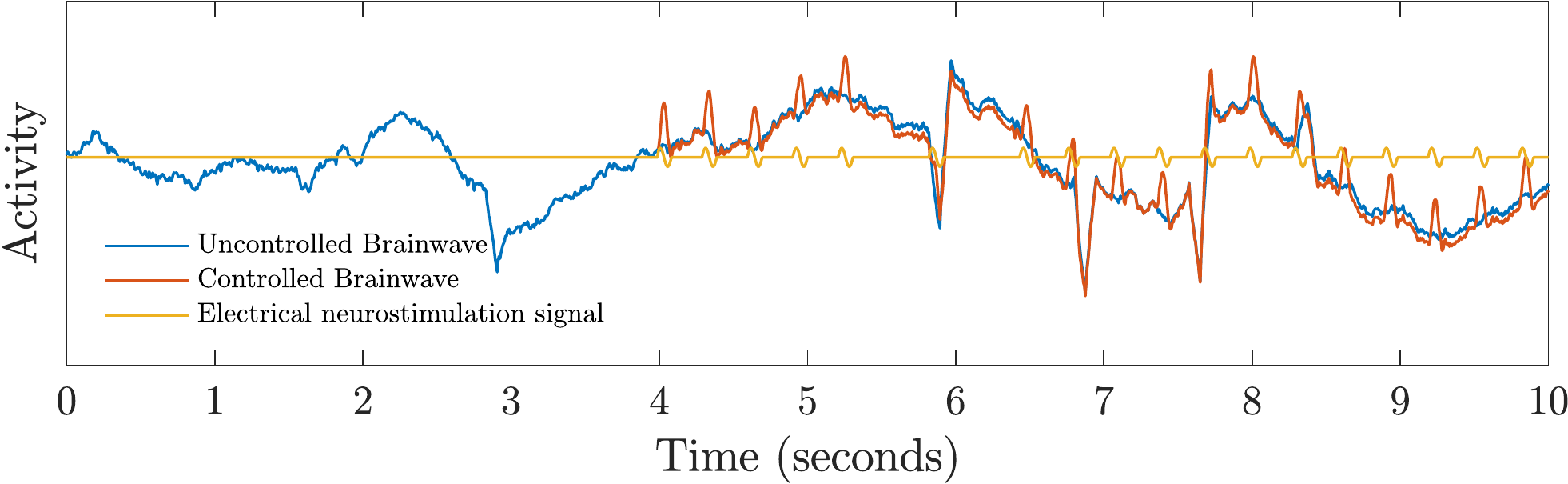}
    \caption{Simulated seizure is depicted (in blue) versus an attempted controlled signal (in red) under an event-triggered open-loop stimulation strategy (in yellow).}
    \label{fig:rns}
\end{figure}

Whereas in the Experiment 1 we focused on open-loop strategies that start at fully pre-specified instances of time, in practice, neurostimulation devices such as those used in TMS implement the so-called event-triggered open-loop stimulation strategy. More specifically, the (open-loop) stimulation strategy is only activated when a phenomenon of interest is observed that can be described by the previously obtained data. In what follows, we assumed that, when activated due to the detection of seizure-like activity, the same open-loop stimulation strategy as presented in Fig.~\ref{fig:openloopseizure} with twice the amplitude is deployed -- see Fig.~\ref{fig:rns}. In practice, there are a variety of such seizure-like activity detectors that can be considered (e.g., the line length feature that measures the length of the line described by the activity during a window of time~\cite{linesearch}). Notwithstanding, we obtained a similar conclusion as when open-loop strategies were not triggered by an event.

\subsection*{Experiment 3: Closed-Loop Electrical Neurostimulation using MPC on FOS}
Finally, we implemented our proposed TMS strategy using MPC on the FOS models described earlier in this paper. More precisely, for the cost function in \eqref{eq:mpc123}, we utilized \mbox{$Q_k=10I_n$}, $R_k=I_{n_u}$, and $c_k=0_{n_u\times 1}$, to emphasize minimizing the overall energy in the measured brainwaves, while penalizing slightly for overly aggressive stimulation. Furthermore, we included a safety linear constraint of \mbox{$-1\leq u_k\leq 1$}. Our predictive model was based on a $(p=4)-$step (2.5 ms) MVAR predictive model approximation of the FOS plant, with a \mbox{$(P=32)-$}step (20 ms) prediction horizon and \mbox{$(M=8)-$}step (5 ms) control horizon.
\begin{figure}[!ht]
    \centering
    \includegraphics[width=0.5\textwidth]{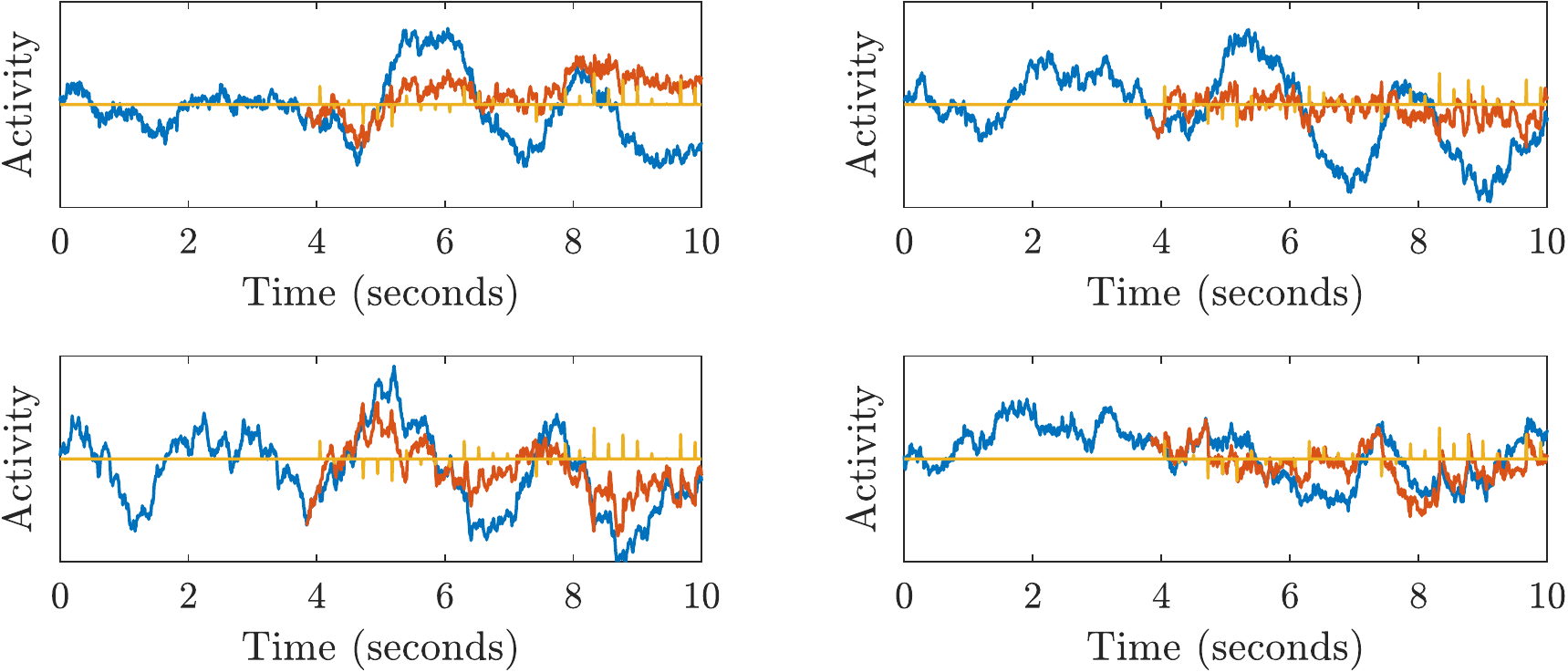}
    \caption{Simulated seizure is depicted (in blue) versus controlled signals (in red) for all channels, under our proposed closed-loop stimulation strategy (in yellow).}
    \label{fig:mpc}
\end{figure}

The results presented in Fig.~\ref{fig:mpc} include all 4 controlled channels controlled simultaneously, since failure to drive any one of them towards a normal range would imply failure in the seizure mitigation objective. We can see that the proposed strategy using TMS achieves the desired goal and implicitly provides us with a detector, given that it tends to provide virtually no stimulation except during very brief periods of time, at which point only (time-varying) impulse-like stimuli are deployed.

\begin{remark}
Any control scheme that requires computationally demanding large-scale optimization methods to be involved, will, in turn, require some form of approximation in order to increase computational efficiency and to enable true real-time control of the system. As such, it is imperative to understand the trade-offs in computational performance involved when using such approximations, which naturally depend from system to system. 
\end{remark}
\begin{remark}
The actual dynamics of the brain are highly nonlinear and time-varying. For this reason, it is crucial to re-identify the fractional-order parameters in our predictive model in a real-life implementation of our proposed strategy. An empirical analysis of parameter tuning via trial and error suggests that an LTI predictive model approximation of the FOS plant with non-augmented state (i.e., an MVAR($p$) predictive model with $p=1$) is inadequate to achieve satisfactory levels of performance in the context of seizure mitigation. However, it also suggests that the memory of the predictive model can be limited to $p=8$ past time steps (ranging 0.5 ms) for the conducted experiments. Anything beyond that led to negligible gains in performance at a high computational cost.
\end{remark}
\begin{remark}
It should also be noted that in the proposed FOS-MPC stimulation strategy, there are still some design parameters that need to be manually calibrated, such as the prediction horizon $P$, the control horizon $M$, the memory horizon $p$, and the input energy penalization weight $\varepsilon\geq 0$. Notwithstanding these considerations, there is a considerable theoretical foundation dedicated to studying the design of MPC algorithms that achieve stability, robustness, and other performance guarantees~\cite{964698,1025358,robustMPC,Pannek}. This body of results may be used to guide and systematize the parameter calibration stage under a sound and justifiable basis.
\end{remark}

\section{Conclusions and Future Work}\label{sec:conclusions}

We presented a methodological framework towards \mbox{real-time} feedback control with constraints for neurophysiological systems. Specifically, we pedagogically introduced a model predictive control (MPC) approach when the neurophysiological process can be modeled by a fractional-order system. In doing so, we focused on TMS for epilepsy, and using systems with seizure-like characteristics, we showed that the stimulation strategies obtained by the proposed framework enabled us to annihilate and mitigate epileptic seizures. Although we have focused mainly on TMS for epilepsy, we believe that the proposed framework can be readily applied to other forms of neurostimulation with an adequate change in the optimization problem (i.e., in the objective function and constraints).

\bibliographystyle{IEEEtran}
\bibliography{IEEEabrv,mybibfile}

\end{document}